\documentclass[a4paper,11pt]{amsart}
\usepackage{graphicx}
\usepackage{mathptmx}
\usepackage{amsmath}
\usepackage{amssymb}
\usepackage{enumitem}
\usepackage{xcolor}
\usepackage{xparse}
\NewDocumentCommand{\eulerian}{omm}
 {%
  \genfrac<>{0pt}{}{#2}{#3}%
  \IfValueT{#1}{_{\!#1}}%
 }

 \newmuskip\pFqmuskip

\newcommand*\pFq[6][8]{%
  \begingroup 
  \pFqmuskip=#1mu\relax
  \mathchardef\normalcomma=\mathcode`,
  \mathcode`\,=\string"8000
  \begingroup\lccode`\~=`\,
  \lowercase{\endgroup\let~}\pFqcomma
  {}_{#2}F_{#3}{\left(\genfrac..{0pt}{}{#4}{#5}\bigg|#6\right)}%
  \endgroup
}
\newcommand{\pFqcomma}{{\normalcomma}\mskip\pFqmuskip}

\newtheorem{theorem}{Theorem}

\begin{document}

\title[On $r$-truncated degenerate Poisson Random Variables]{ Some results on $r$-truncated degenerate Poisson Random Variables}

\author{Taekyun  Kim}
\address{Department of Mathematics, Kwangwoon University, Seoul 139-701, Republic of Korea}
\email{tkkim@kw.ac.kr}

\author{DAE SAN KIM}
\address{Department of Mathematics, Sogang University, Seoul 121-742, Republic of Korea}
\email{dskim@sogang.ac.kr*}

\author{Si-Hyeon Lee}
\address{Department of Mathematics, Kwangwoon University, Seoul 139-701, Republic of Korea}
\email{ugug11@naver.com}

\author{Seong-Ho Park}
\address{Department of Mathematics, Kwangwoon University, Seoul 139-701, Republic of Korea}
\email{abcd2938471@kw.ac.kr}

\author{Lee-Chae Jang }
\address{Graduate School of Education, Konkuk University, Seoul 143-701, Republic of Korea}
\email{lcjang@konkuk.ac.kr}


\maketitle

\begin{abstract}
The zero-truncated Poisson distributions are certain discrete probability distributions whose supports are the set of positive integers, which are also known as the conditional Poisson distributions or the positive Poisson distributions.
Recently, as a natural extension of those distributions, Kim-Kim studied the zero-truncated degenerate Poisson distributions. In this paper, we introduce the $r$-truncated degenerate Poisson random variable with parameter $\alpha > 0$, whose probability mass function is given by $p_{\lambda,r}(i)=\cfrac{(1)_{i,\lambda}}{e_{\lambda}(\alpha)-e_{\lambda,r}(\alpha)}\frac{\alpha^{i}}{i!},\quad(i=r+1,r+2,r+3,\dots)$, and investigate various properties of this random variable.
\end{abstract}

\section{Introduction}
It is well known that a random variable $X$, taking on one of the values $0,1,2,\dots$, is said to be the Poisson random variable with parameter $\alpha>0$, if the probability mass function of $X$ is given by
\begin{displaymath}
p(i)=P\{X=i\}=e^{-\alpha}\frac{\alpha^{i}}{i!},\quad(i=0,1,2,\dots),\quad(\mathrm{see}\ [11,13]).
\end{displaymath}
We note that a Poisson random variable indicates how many events occured in a given period of time.
Further, a random variable $X_{r}$, taking on one of the values $r+1,r+2,\dots$, is called the  $r$-truncated Poisson random variable with parameter $\alpha >0$, if the probability mass function of $X_{r}$ is given by
\begin{displaymath}
p_{r}(i)=P\{X_{r}=i\}=\frac{1}{e^{\alpha}-\sum_{l=0}^{r}\frac{\alpha^l}{l!}}\frac{\alpha^{i}}{i!},
\end{displaymath}
where $i=r+1,r+2,r+3,\dots$, for $r \ge 0$. In the special case of $r=0$, we obtain the zero-truncated  Poisson
distribution.
\par
\indent On the other hand, the degenerate Poisson random variables with parameter $\alpha>0$, whose probability mass function is given by (see (1))
\begin{displaymath}
p_{\lambda}(i)=e_{\lambda}^{-1}(\alpha)\frac{\alpha^{i}}{i!}(1)_{i,\lambda}\quad(i=0,1,2,\dots),
\end{displaymath}
were studied by Kim-Kim-Jang-Kim in [10].
In addition, the zero-truncated degenerate Poisson random variables with parameter $\alpha>0$, whose probability mass function is given by
\begin{displaymath}
p_{\lambda,0}(i)=\frac{1}{e_{\lambda}(\alpha)-1}\frac{\alpha^{i}}{i!}(1)_{i,\lambda}\quad(i=1,2,3,\dots),
\end{displaymath}
were studied in [9]. \par
The aim of this paper is to generalize the results on the zero-truncated degenerate Poisson distribution in [9]  to the case of the $r$-truncated degenerate Poisson distribution as a natural extension of the $r$-truncated Poission distribution. We will obtain, among other things, its expectation, its variance, its $n$-th moment, and its cumulative distribution function. \par
One motivation for this research is its potential applications to the cornavirus pandemic. It has been spreading unpredictably around the world, terrorizing many people. Although several vaccines for the coronavirus have been developed and many people are getting vaccinated, they have many unexpected side effects as well.
We would like to predict stability of the coronavirus vaccines after $r$-days the vaccines were shot. For this purpose, we study the $r$-truncated degenerate Poisson random variables (see (6)), which has the `degenerate factor' $\lambda$ reflecting abnormal situations.
Indeed, we think that Theorem 5 is useful in predicting the probability of the coronavirus vaccines becoming stable after the $r$-days of getting vaccinated. Another motivation is applications of various probabilistic methods in studying special numbers and polynomials arising from combinatorics and number theory. For this, we let the reader refer to [5-10,14 and the references therein]. For the rest of this section, we recall the necessary facts that are needed throughout this paper. \par
For any $\lambda\in\mathbb{R}$, the degenerate exponential functions are defined by
\begin{equation}
e_{\lambda}^{x}(t)=(1+\lambda t)^{\frac{x}{\lambda}}=\sum_{k=0}^{\infty}(x)_{k,\lambda}\frac{t^{k}}{k!},\quad(\mathrm{see}\ [1-3,5-10,14]),\label{1}
\end{equation}
where $(x)_{0,\lambda}=1$, $(x)_{n,\lambda}=x(x-\lambda)\cdots(x-(n-1)\lambda)$, $(n\ge 1)$. \par
\noindent In particular, for $x=1$, we denote $e_{\lambda}^{1}(t)$ by $e_{\lambda}(t)$. \par
\noindent In addition, it is convenient to introduce $e_{\lambda,r}(t)$, for any integer $r \ge 0$, given by
\begin{equation}
e_{\lambda,r}(t)=\sum_{k=0}^{r}(1)_{k,\lambda}\frac{t^{k}}{k!}. \label{2}	
\end{equation} \par
The Bell polynomials are given by
\begin{equation}
e^{x(e^{t}-1)}=\sum_{n=0}^{\infty}\mathrm{Bel}_{n}(x)\frac{t^{n}}{n!},\quad(\mathrm{see}\ [4,12]).\label{3}
\end{equation}
In light of \eqref{3}, the degenerate Bell polynomials are defined in [9] by
\begin{equation}
e_{\lambda}^{-1}(x)e_{\lambda}(xe^{t})=\sum_{n=0}^{\infty}\mathrm{Bel}_{n,\lambda}(x)\frac{t^{n}}{n!}. \label{4}
\end{equation}
In particular, for $x=1$, $\mathrm{Bel}_{n,\lambda}=\mathrm{Bel}_{n,\lambda}(1)$, $(n\ge 0)$, are called the degenerate Bell numbers. Note that
\begin{displaymath}
	\lim_{\lambda\rightarrow 1}\mathrm{Bel}_{n,\lambda}(x)=\mathrm{Bel}_{n}(x),\quad(n\ge 0).
\end{displaymath} \par
Let $X$ be a Poisson random variable with parameter $\alpha>0$. Then the moment of $X$ is given by
\begin{displaymath}
	E[X^{n}]=\mathrm{Bel}_{n}(\alpha),\quad(n\ge 0),\quad(\mathrm{see}\ [11,13]).
	\end{displaymath}
	Let $Y$ be a discrete random variable taking values in the nonnegative integers. Then the probability generating function of $Y$ is given by
	\begin{displaymath}
	F(t)=E[t^{Y}]=\sum_{i=0}^{\infty}p(i)t^{i},\quad(\mathrm{see}\ [11,13]),
	\end{displaymath}
	where $p(i)=P\{Y=i\}$ is the probability mass function of $Y$. \par
The degenerate Stirling numbers of the second kind $S_{2,\lambda}(n,k)$ are defined by
\begin{displaymath}
(x)_{n,\lambda}=\sum_{l=0}^{n}S_{2,\lambda}(n.l)(x)_{l},\quad(n\ge 0),\quad(\mathrm{see}\ [5]),
\end{displaymath}
or equivalently by
\begin{displaymath}
\frac{1}{k!}\bigg(e_{\lambda}(t)-1\bigg)^{k}=\sum_{n=k}^{\infty}S_{2,\lambda}(n,k)\frac{t^{n}}{n!}.
\end{displaymath}
Here $(x)_{0}=1,\ (x)_{n}=x(x-1)(x-2)\cdots(x-n+1)$, $(n\ge 1)$. \par
Now, we consider the {\it{$r$-truncated degenerate Stirling numbers of the second kind}} $S_{2,\lambda}^{[r]}(n,kr+k)$ given by
\begin{equation}
\frac{1}{k!}\bigg(e_{\lambda}(t)-e_{\lambda,r}(t)\bigg)^{k}=\sum_{n=kr+k}^{\infty}S_{2,\lambda}^{[r]}(n,kr+k)\frac{t^{n}}{n!}.\label{5}	
\end{equation}
Note that
\begin{displaymath}
S_{2,\lambda}^{[0]}(n,k)=S_{2,\lambda}(n,k),\quad(\mathrm{see}\ [5]).
\end{displaymath}

\section{$r$-truncated degenerate Poisson random variables}
	
The random variable $X_{\lambda,r}$ is called the {\it{$r$-truncated degenerate Poisson random variable with parameter $\alpha>0$}}, if the probability mass function of $X_{\lambda,r}$ is given by
\begin{equation}
p_{\lambda,r}(k)=\cfrac{(1)_{k,\lambda}}{e_{\lambda}(\alpha)-e_{\lambda,r}(\alpha)}\frac{\alpha^{k}}{k!}, \label{6}
\end{equation}
where $k=r+1,r+2,\dots$, with $r\ge 0$. \par
Here we must observe that
\begin{align*}
\sum_{k=r+1}^{\infty}p_{\lambda,r}(k) &= \frac{1}{e_{\lambda}(\alpha)-e_{\lambda,r}(\alpha)}\sum_{k=r+1}^{\infty}\frac{\alpha^k}{k!}(1)_{k,\lambda} \\
&= \frac{1}{e_{\lambda}(\alpha)-e_{\lambda,r}(\alpha)}(e_{\lambda}(\alpha)-e_{\lambda,r}(\alpha))\\
&=1.
\end{align*}
In addition, we note that
\begin{displaymath}
\lim_{\lambda\rightarrow 0}p_{\lambda,r}(k)=\frac{1}{e^{\alpha}-\sum_{l=0}^{r}\frac{\alpha^{l}}{l!}}\frac{\alpha^{k}}{k!}
\end{displaymath}
is the probability mass function of the $r$-truncated Poisson random variable with parameter $\alpha>0$. \par
Let us assume that $X_{\lambda,r}$ is the $r$-truncated degenerate Poisson random variable with parameter $\alpha>0$. Then we note that the expectation of $X_{\lambda,r}$ is given by
\begin{align}
E[X_{\lambda,r}]\ &=\ \sum_{n=r+1}^{\infty}np_{\lambda,r}(n)=\frac{1}{e_{\lambda}(\alpha)-e_{\lambda,r}(\alpha)}\sum_{n=r+1}^{\infty}n\frac{\alpha^{n}}{n!}(1)_{n,\lambda}\label{7}\\
&=\ \frac{\alpha} {e_{\lambda}(\alpha)-e_{\lambda,r}(\alpha)}\sum_{n=r}^{\infty}\frac{\alpha^{n}}{n!}(1)_{n+1,\lambda} \nonumber	\\
&=\ \frac{\alpha}{e_{\lambda}(\alpha)-e_{\lambda,r}(\alpha)}\bigg(\sum_{n=r+1}^{\infty}\frac{\alpha^{n}}{n!}(1)_{n+1,\lambda}+\frac{\alpha^{r}}{r!}(1)_{r+1,\lambda}\bigg) \nonumber \\
&=\ \frac{\alpha}{e_{\lambda}(\alpha)-e_{\lambda,r}(\alpha)}\bigg(\sum_{n=r+1}^{\infty}\frac{\alpha^{n}}{n!}(1)_{n,\lambda}(1-n\lambda)+\frac{\alpha^{r}}{r!}(1)_{r+1,\lambda}\bigg) \nonumber\\
&=\ \alpha-\frac{\alpha\lambda}{e_{\lambda}(\alpha)-e_{\lambda,r}(\alpha)}\sum_{n=r+1}^{\infty}n\frac{\alpha^{n}}{n!}(1)_{n,\lambda}+ \frac{\alpha}{e_{\lambda}(\alpha)-e_{\lambda,r}(\alpha)}\frac{\alpha^{r}}{r!}(1)_{r+1,\lambda} \nonumber \\
&=\ \alpha-\alpha\lambda E[X_{\lambda,r}]+ \frac{\alpha}{e_{\lambda}(\alpha)-e_{\lambda,r}(\alpha)}\frac{\alpha^{r}}{r!}(1)_{r+1,\lambda}\nonumber.
\end{align}
From \eqref{7}, we obtain the following theorem.
\begin{theorem}
Let $X_{\lambda,r}$ be the $r$-truncated degenerate Poisson random with parameter $\alpha>0$. Then, for $\lambda\ne-\frac{1}{\alpha}$, the expectation of $X_{\lambda,r}$ is given by
\begin{displaymath}
E[X_{\lambda,r}]=\frac{1}{1+\alpha\lambda}\bigg[\alpha+\frac{\alpha^{r+1}}{r!} \frac{(1)_{r+1,\lambda}}{e_{\lambda}(\alpha)-e_{\lambda,r}(\alpha)}\bigg].
\end{displaymath}
\end{theorem}
Now, we observe that
\begin{align}
E[X_{\lambda,r}^{2}]\ &=\ \sum_{n=k+1}^{\infty}n^{2}p_{\lambda,r}(n)\ =\ \frac{1}{e_{\lambda}(\alpha)-e_{\lambda,r}(\alpha)}\sum_{n=r+1}^{\infty}n^{2}\frac{\alpha^{n}}{n!}(1)_{n,\lambda} \label{8} \\
&=\ \frac{1}{e_{\lambda}(\alpha)-e_{\lambda,r}(\alpha)}\bigg(\sum_{n=r+1}^{\infty}\frac{\alpha^{n}}{(n-2)!}(1)_{n,\lambda}+\sum_{n=r+1}^{\infty}\frac{\alpha^{n}}{(n-1)!}(1)_{n,\lambda}\bigg)\nonumber \\
&=\ \frac{1}{e_{\lambda}(\alpha)-e_{\lambda,r}(\alpha)}\bigg(\sum_{n=r+2}^{\infty}\frac{\alpha^{n}}{(n-2)!}(1)_{n,\lambda}+\sum_{n=r+1}^{\infty}\frac{\alpha^{n}}{(n-1)!}(1)_{n,\lambda}\bigg) \nonumber \\
&\quad +\frac{1}{e_{\lambda}(\alpha)-e_{\lambda,r}(\alpha)}\frac{\alpha^{r+1}}{(r-1)!}(1)_{r+1,\lambda} \nonumber \\
&=\ \frac{1}{e_{\lambda}(\alpha)-e_{\lambda,r}(\alpha)}\bigg(\alpha\sum_{n=r+1}^{\infty}\frac{\alpha^{n}}{(n-1)!}(1)_{n+1,\lambda}+\sum_{n=r+1}^{\infty}\frac{n\alpha^{n}}{n!}(1)_{n,\lambda}\bigg) \nonumber \\
&\quad +\frac{1}{e_{\lambda}(\alpha)-e_{\lambda,r}(\alpha)}r\frac{\alpha^{r+1}}{r!}(1)_{r+1,\lambda} \nonumber \\
&=\ \frac{\alpha}{e_{\lambda}(\alpha)-e_{\lambda,r}(\alpha)}\sum_{n=r+1}^{\infty}\frac{n\alpha^{n}}{n!}(1)_{n,\lambda}(1-n\lambda)+E[X_{\lambda,r}]\nonumber \\
&\quad +\frac{r\alpha^{r+1}(1)_{r+1,\lambda}}{r!\big(e_{\lambda}(\alpha)-e_{\lambda,r}(\alpha)\big)} \nonumber \\
&=\ \alpha E[X_{\lambda,r}]-\lambda\alpha E[X_{\lambda,r}^{2}]+E[X_{\lambda,r}]+\frac{\alpha^{r+1}}{r!}\frac{r(1)_{r+1,\lambda}}{e_{\lambda}(\alpha)-e_{\lambda,r}(\alpha)}.\nonumber
\end{align}
From \eqref{8}, we note that
\begin{equation}
E[X_{\lambda,r}^{2}]=\frac{\alpha+1}{1+\lambda\alpha}E[X_{\lambda,r}]+\frac{1}{1+\lambda\alpha} \frac{r(1)_{r+1,\lambda}}{e_{\lambda}(t)-e_{\lambda,r}(\alpha)}\frac{\alpha^{r+1}}{r!},\label{9}	
\end{equation}
where $\lambda\ne-\frac{1}{\alpha}$. \par
For $\lambda\ne-\frac{1}{\alpha}$, we note that the variance of $X_{\lambda,r}$ is given by
\begin{align}
&\mathrm{Var}(X_{\lambda,r})\ =\ E[X_{\lambda,r}^{2}]-\big(E[X_{\lambda,r}]\big)^{2} \label{10}\\
&=\ E[X_{\lambda,r}]\bigg(\frac{\alpha+1}{1+\lambda\alpha}-E[X_{\lambda,r}]\bigg)+ \frac{1}{1+\lambda\alpha}\frac{r(1)_{r+1,\lambda}}{e_{\lambda}(\alpha)-e_{\lambda,r}(\alpha)}\frac{\alpha^{r+1}}{r!} \nonumber \\
&=\ \frac{1}{1+\alpha\lambda}\bigg(\alpha+\frac{\alpha^{r+1}}{r!}\frac{(1)_{r+1,\lambda}}{e_{\lambda}(\alpha)-e_{\lambda,r}(\alpha)}\bigg)\bigg(\frac{\alpha+1}{1+\alpha\lambda}-\frac{1}{1+\alpha\lambda}\bigg(\alpha+\frac{\alpha^{r+1}}{r!}\frac{(1)_{r+1,\lambda}}{e_{\lambda}(\alpha)-e_{\lambda,r}(\alpha)}\bigg)\bigg)\nonumber\\
&\quad +\frac{1}{1+\alpha\lambda}\frac{r(1)_{r+1,\lambda}}{e_{\lambda}(\alpha)-e_{\lambda,r}(\alpha)}\frac{\alpha^{r+1}}{r!} \nonumber \\
&=\ \frac{1}{(1+\alpha\lambda)^{2}}\bigg(\alpha+\frac{\alpha^{r+1}}{r!}\frac{(1)_{r+1,\lambda}}{e_{\lambda}(\alpha)-e_{\lambda,r}(\alpha)}\bigg)\bigg(1-\frac{\alpha^{r+1}}{r!}\frac{(1)_{r+1,\lambda}}{e_{\lambda}(\alpha)-e_{\lambda,r}(\alpha)}\bigg)\nonumber \\
&\quad+\frac{1}{1+\alpha\lambda}\frac{r(1)_{r+1,\lambda}}{e_{\lambda}(\alpha)-e_{\lambda,r}(\alpha)}\frac{\alpha^{r+1}}{r!}.\nonumber
\end{align}
Therefore, by \eqref{10}, we obtain the following theorem.
\begin{theorem}
Let $X_{\lambda,r}$ be the $r$-truncated degenerate Poisson random variable with parameter $\alpha>0$. For $\lambda\ne-\frac{1}{\alpha}$, we have
\begin{align*}
\mathrm{Var}(X_{\lambda,r})&=\frac{1}{(1+\alpha\lambda)^{2}}\bigg(\alpha+\frac{\alpha^{r+1}}{r!}\frac{(1)_{r+1,\lambda}}{e_{\lambda}(\alpha)-e_{\lambda,r}(\alpha)}\bigg)\bigg(1-\frac{\alpha^{r+1}}{r!}\frac{(1)_{r+1,\lambda}}{e_{\lambda}(\alpha)-e_{\lambda,r}(\alpha)}\bigg) \\
&+\frac{1}{1+\alpha\lambda}\frac{r(1)_{r+1,\lambda}}{e_{\lambda}(\alpha)-e_{\lambda,r}(\alpha)}\frac{\alpha^{r+1}}{r!}.
\end{align*}
\end{theorem}
Let us consider the generating function of the moments of the $r$-truncated degenerate Poisson random variable with parameter $\alpha>0$. Then, from \eqref{4}, we have
\begin{align}
\sum_{n=0}^{\infty}E[X_{\lambda,r}^{n}]\frac{t^{n}}{n!}\ &=\ E[e^{X_{\lambda,r}t}]\ =\ \sum_{m=r+1}^{\infty}e^{mt}p_{\lambda,r}(m) \label{11} \\
&=\ \frac{1}{e_{\lambda}(\alpha)-e_{\lambda,r}(\alpha)}\sum_{m=r+1}^{\infty}\frac{\alpha^{m}}{m!}(1)_{m,\lambda}e^{mt}\nonumber\\
&=\ \frac{1}{e_{\lambda}(\alpha)-e_{\lambda,r}(\alpha)}\bigg(e_{\lambda}(\alpha e^{t})-\sum_{m=0}^{r}\frac{\alpha^{m}}{m!}(1)_{m,\lambda}e^{mt}\bigg) \nonumber	\\
&=\ \frac{e_{\lambda}(\alpha)}{e_{\lambda}(\alpha)-e_{\lambda,r}(\alpha)} \bigg(e_{\lambda}^{-1}(\alpha)e_{\lambda}(\alpha e^{t})-e_{\lambda}^{-1}(\alpha)\sum_{m=0}^{r}\frac{\alpha^{m}}{m!}(1)_{m,\lambda}e^{mt}\bigg)\nonumber\\
&=\ \frac{e_{\lambda}(\alpha)}{e_{\lambda}(\alpha)-e_{\lambda,r}(\alpha)}\bigg(\sum_{n=0}^{\infty}\mathrm{Bel}_{n,\lambda}(\alpha)\frac{t^{n}}{n!}-\sum_{n=0}^{\infty}e_{\lambda}^{-1}(\alpha)\sum_{m=0}^{r}\frac{\alpha^{m}}{m!}(1)_{m,\lambda}m^{n}\frac{t^{n}}{n!}\bigg) \nonumber \\
&=\ \sum_{n=0}^{\infty}\frac{e_{\lambda}(\alpha)}{e_{\lambda}(\alpha)-e_{\lambda,r}(\alpha)}\bigg(\mathrm{Bel}_{n,\lambda}(\alpha)-e_{\lambda}^{-1}(\alpha)\sum_{m=0}^{r}\frac{\alpha^{m}}{m!}(1)_{m,\lambda}m^{n}\bigg)\frac{t^{n}}{n!}.\nonumber
\end{align}
Thus, by \eqref{11}, we get the next theorem.
\begin{theorem}
Let $X_{\lambda,r}$ be the $r$-truncated degenerate Poisson random variable with parameter $\alpha>0$. For $n\ge0$, we have
\begin{displaymath}
E[X_{\lambda,r}^{n}]= \frac{e_{\lambda}(\alpha)}{e_{\lambda}(\alpha)-e_{\lambda,r}(\alpha)}\bigg(\mathrm{Bel}_{n,\lambda}(\alpha)-e_{\lambda}^{-1}(\alpha)\sum_{m=0}^{r}\frac{\alpha^{m}}{m!}(1)_{m,\lambda}m^{n}\bigg).
\end{displaymath}
\end{theorem}
For $x\ge r+1$, we note that the cumulative distribution function is given by
\begin{align}
F_{X_{\lambda,r}}(x)\ &=\ P\{X_{\lambda,r}\le x\}\ =\ \sum_{k=r+1}^{[x]}p_{\lambda,r}(k) \label{12}\\
&=\ \frac{1}{e_{\lambda}(\alpha)-e_{\lambda,r}(\alpha)}\sum_{k=r+1}^{[x]}\frac{\alpha^{k}}{k!}(1)_{k,\lambda} \nonumber \\
&=\ \frac{1}{e_{\lambda}(\alpha)-e_{\lambda,r}(\alpha)}\bigg(\sum_{k=0}^{[x]}\frac{\alpha^{k}}{k!}(1)_{k,\lambda}-\sum_{k=0}^{r}\frac{\alpha^{k}}{k!}(1)_{k,\lambda}\bigg). \nonumber
\end{align}

From \eqref{12} and \eqref{2}, we can derive the following equation.
\begin{equation}
F_{X_{\lambda,r}}(x)=\frac{1}{e_{\lambda}(\alpha)-e_{\lambda,r}(\alpha)}\big(e_{\lambda,[x]}(\alpha)-e_{\lambda,r}(\alpha)\big). \label{13}
\end{equation}
\begin{theorem}
Assume that $X_{\lambda,r}$ is the $r$-truncated degenerate Poisson random variable with parameter $\alpha>0$. For $x\ge r+1$, the cumulative distribution function of $X_{\lambda,r}$ is given by
\begin{displaymath}
F_{X_{\lambda,r}}(x)=\frac{1}{e_{\lambda}(\alpha)-e_{\lambda,r}(\alpha)}\big(e_{\lambda,[x]}(\alpha)-e_{\lambda,r}(\alpha)\big).
\end{displaymath}
\end{theorem}
Let us assume that $X_{\lambda,r}^{(1)},X_{\lambda,r}^{(2)},\dots,X_{\lambda,r}^{(k)}$ are identically independent $r$-truncated degenerate Poisson random variables with parameter $\alpha>0$, and let
\begin{displaymath}
X_{\lambda,r}=\sum_{i=1}^{k}X_{\lambda,r}^{(i)},\quad(k\in\mathbb{N}).
\end{displaymath}
From the probability generating function of random variable, we note that
\begin{align}
E[t^{X_{\lambda,r}^{(i)}}]\ &=\ \sum_{n=r+1}^{\infty}t^{n}P[X_{\lambda,r}^{(i)}=n] \label{14} \\
&=\ \frac{1}{e_{\lambda}(\alpha)-e_{\lambda,r}(\alpha)}\sum_{n=r+1}^{\infty}\frac{\alpha^{n}}{n!}(1)_{n,\lambda}t^{n} \nonumber \\
&=\ \frac{e_{\lambda}(\alpha t)-e_{\lambda,r}(\alpha t)}{e_{\lambda}(\alpha)-e_{\lambda,r}(\alpha)}.\nonumber
\end{align}
By \eqref{14}, we get
\begin{align}
E[t^{X_{\lambda,r}}]\ &=\ \prod_{i=1}^{k}E[t^{X_{\lambda,r}^{(i)}}] \label{15} \\
&=\ \bigg(\frac{1}{e_{\lambda}(\alpha)-e_{\lambda,r}(\alpha)}\bigg)^{k}k!\frac{1}{k!}\big(e_{\lambda}(\alpha t)-e_{\lambda,r}(\alpha t)\big)^{k} \nonumber \\
&=\ \frac{k!}{(e_{\lambda}(\alpha)-e_{\lambda,r}(\alpha))^{k}}\sum_{n=kr+k}^{\infty}S_{2,\lambda}^{[r]}(n,kr+k)\frac{\alpha^{n}t^{n}}{n!}. \nonumber
\end{align}
On the other hand,
\begin{align}
E[t^{X_{\lambda,r}}]\ &=\ E[t^{X_{\lambda,r}^{(1)}+X_{\lambda,r}^{(2)}+\cdots+X_{\lambda,r}^{(k)}}]\label{16} \\
&=\ \sum_{n=kr+r}^{\infty}t^{n}P\bigg[\sum_{i=1}^{k}X_{\lambda,r}^{(i)}=n	\bigg].\nonumber
\end{align}
Therefore, by \eqref{15} and \eqref{16}, we obtain the following theorem.
\begin{theorem}
Let $X_{\lambda,r}^{(1)},X_{\lambda,r}^{(2)},\dots,X_{\lambda,r}^{(k)}$ be identically independent $r$-truncated degenerate Poisson random variables with parameter $\alpha>0$, and let $\displaystyle X_{\lambda,r}=\sum_{i=1}^{k}X_{\lambda,r}^{(i)}\displaystyle$. Then the probability for $X_{\lambda,r}$ is given by
\begin{displaymath}
P[X_{\lambda,r}=n]=\left\{\begin{array}{ccc}
\cfrac{k!}{(e_{\lambda}(\alpha)-e_{\lambda,r}(\alpha))^{k}}\frac{\alpha^{n}}{n!}S_{2,\lambda}^{[r]}(n,kr+k), & \textrm{if $n\ge kr+r$}, \\
0, & \textrm{otherwise.}
\end{array}\right.
\end{displaymath}
\end{theorem}

\section{Conclusion}

We generalized the results on the zero-truncated degenerate Poisson distribution to the case of the $r$-truncated degenerate Poisson distribution, with its potential applications to the cornavirus pandemic and applications of probabilistic methods to the study of some special numbers and polynomials in mind. \par
Let $X_{\lambda,r}$ be the $r$-truncated degenerate Poisson random variable with parameter $\alpha$.
Then, for the random variable $X_{\lambda,r}$, we derived its expectation, its variance, its $n$-th moment, and its cumulative distribution function. In addition, we obtained two different expressions for the probability generating function of a finite sum of independent $r$-truncated degenerate Poisson random variables with equal parameters. \\
\indent As one of our future projects, we would like to continue this line of research, namely to explore applications of various methods of probability theory to science, engineering and social science, and to the study of some special polynomials and numbers.

\vspace{0.1in}
{\bf{Acknowledgements:}}  Not applicable.
\vspace{0.1in}

{\bf{Funding:}} Not applicable.
\vspace{0.1in}

{\bf {Availability of data and material:}}
Not applicable.

\vspace{0.1in}

{\bf {Competing interests:}}
The authors declare no conflict of interest.

\vspace{0.1in}

{\bf{Author Contributions:}} D.S.K. and T.K. wrote the paper; L.-C.J. and S.-H. L. and S.-H. P. checked the results of the paper and typed the paper. All authors have read and agreed to the published version of the manuscript.

\end{document}